\input amstex
\documentstyle{amsppt}
\document
\magnification=1200
\NoBlackBoxes
\nologo
\pageheight{18cm}

\centerline{\bf MIRROR SYMMETRY AND QUANTIZATION}

\smallskip

\centerline{\bf OF ABELIAN VARIETIES}

\medskip

\centerline{\bf Yu.~I.~Manin}

\medskip

\centerline{\it Max--Planck--Institut f\"ur Mathematik, Bonn}

\bigskip

\centerline{\bf 0. Introduction}

\medskip

{\bf 0. Plan of the paper.} This paper consists of two sections discussing
various aspects of commutative and non--commutative geometry
of tori and abelian varieties. 

\smallskip

In the first section,
we present a new definition of mirror symmetry
for abelian varieties and, more generally,
complex and $p$--adic tori, that is, spaces
of the form $T/B$ where $T$ is an algebraic group isomorphic
to a product of multiplicative groups, $K$ is
a complete normed field, and $B\subset T(K)$ is
a discrete subgroup of maximal rank in it. We also
check its compatibility with other definitions
discussed in the literature.

\smallskip

In the second section, we develop an approach
to the quantization of abelian varieties first
introduced in [Ma1], namely, via theta
functions on {\it non--commutative, or quantum, tori} endowed
with a discrete period lattice. These
theta functions satisfy a functional equation which
is a generalization of the classical one, in particular,
involve a multiplier. Since multipliers cease
to be central in the quantum case, one must decide
where to put them. In [Ma1] only one--sided multipliers
were considered. As a result, the product of two theta functions
in general was not a theta function. Here we suggest a partial
remedy to this problem by introducing two--sided multipliers.
The resulting space of theta functions 
possesses partial multiplication and has sufficiently
rich functorial properties so that rudiments of
Mumford's theory ([Mu]) can be developed. Main results
of [Ma1] are reproduced here in a generalized form,
so that this paper can be read independently.

\smallskip

We will now briefly describe a broader picture into which
this work fits.

\medskip

{\bf 1. On mirror symmetry.} In this paper mirror symmetry is 
understood as a binary
relation between (weak) Calabi--Yau
manifolds endowed with some additional data.
A projective (or compact complex) manifold $V$ is a weak Calabi--Yau,
if it admits nowhere vanishing global volume form.
Additional data which are commonly
considered are of two types.

\smallskip

A) A symplectic or complexified symplectic structure $\omega_V$
on $V$, which is ``sufficiently large'' in the case of complex base field.

\smallskip

B) A cusp $c_V$ in the moduli space (or rather stack) of 
deformations of $V$, that is,
a neighborhood of a point of ``maximal degeneration'' to which
$V$ belongs.

\smallskip

The mirror partnership relation between $(V,c_V,\omega_V)$
and $(W,c_W,\omega_W)$ consists in a host of identifications
(partly conjectural)
of various structures that can be produced starting with such
triples. In the case of strict Calabi--Yau's, this includes
an identification of two Frobenius manifolds:
quantum cohomology of $(V,\omega_V)$ and a germ of the extended moduli
space of $W$ with its flat structure determined by $c_W$,
and similarly with roles of $V$ and $W$ reversed.
Generally, one expects also a representation of $V$ and $W$
as dual real Lagrangian torus fibrations over a common base,
with a rich structure of Fourier--Mukai transform
connecting Lagrangian/complex analytic objects on both sides.
For more details on this, see original works
[MirS1], [MirS2], [StYZ], [Ko], [Giv], [Bar], [LYZ], and a report 
[Ma3].

\smallskip

Most important testing ground
for all levels of mirror correspondence is furnished
by toric mirrors introduced and studied by V.~Batyrev.
In the simplest version, this construction
looks as follows (see [Bat]).

\smallskip

Let $T$ be an $n$--dimensional algebraic torus, that is an algebraic group
which is isomorphic to a product of $n$ multiplicative groups $\bold{G}_m$.
It determines (and is functorially determined by any of) the two
free abelian groups of rank $n$: its character group
$M_T:=\roman{Hom}\,(T,\bold{G}_m)$ and its group
of one--parametric subgroups $N_T:=\roman{Hom}\,(\bold{G}_m,T)$.
These groups are naturally dual to each other. Denote by
$T^t$ the dual torus, whose character group is $M_{T^t}=N_T$ and 
respectively $N_{T^t}=M_T$. In the main text of paper,
one of these groups is denoted $H$ and another $H^t.$

\smallskip

Among various toric compactifications of $T$ we are interested in
those for which their anti\-canonical system is ample.
Anticanonical divisors on them are $n$--dimensi\-onal Calabi--Yau
manifolds.

\smallskip

Such toric compactifications $\overline{T}_{\Delta}$ 
are naturally indexed by {\it reflexive polyhedra} 
$\Delta$ in $N_T\otimes\bold{R}$. Via standard convex duality, each such polyhedron determines
the dual polyhedron $\Delta^t$ in $N_{T^t}\otimes\bold{R}$
which is also reflexive. According to Batyrev,
families of anticanonical hypersurfaces in the respective
toric compactifications are expected to be mirror partners:
$$
\overline{T}_{\Delta}\supset \left| -K_{\overline{T}_{\Delta}}\right|\
\Longleftrightarrow\
\left| -K_{\overline{T}^t_{\Delta^t}}\right|\subset
\overline{T}^t_{\Delta^t} .
\eqno(0.1)
$$

The relevant maximally
degenerate CY's are simply divisors at infinity in these
toric compactifications. 

\smallskip

Supplementing $(T,\Delta )$ with additional combinatorial
structure, one can generalize this picture to some Calabi--Yau
complete intersections. 

\smallskip

The main goal of the first section of this paper
is to provide a new definition of mirror symmetry
for abelian varieties (and more generally, complex and $p$--adic tori)
similar to (0.1). We use a ``multiplicative uniformization''
which goes back to Jacobi and which represents $\Cal{A}$
as a quotient of $n$--dimensional algebraic torus $T$
 by a multiplicative
discrete lattice $B\subset T(K).$ 
Consider now two dual algebraic tori endowed with period lattices 
which are explicitly identified,
that is a diagram of the form:
$$
(i,i^t):\ T\leftarrow B\rightarrow T^t\,.
\eqno(0.2)
$$
We will say that pairs $(\Cal{A}:=T(K)/i(B), i^t)$ and
$(\Cal{B}:=T^t(K)/i^t(B), i)$ are mirror dual to each other.
Over complex field, we will relate
$i^t$ (resp. $i$) to a structure similar to $\omega_V$ (resp. $\omega_W$)
above and compare our mirror relation with that of [Gr1], [Gr2], [AP] and [GolLO]. In particular, we will see how this diagram gives rise
to two dual fibrations of the relevant abelian varieties
(or complex tori) by mutually dual Lagrangian real tori over the
same base.

\smallskip

The choice of multiplicative uniformization is not unique,
and it provides the environment for a partial compactification
of the relevant moduli space and  choice  of a maximally
degenerate point at the boundary of moduli space to which $\Cal{A}$
is close, that is, of the relevant cusp
of the moduli space. Roughly speaking, at the boundary some generators of the period lattice are forced to vanish so that the rank of the image of $B$
drops.
To make it more precise,
we can choose a fan $\Phi$ in $N_T\otimes\bold{R}$, construct the dual fan
$\Phi^t$ and study the moduli space of the refined diagrams
$$
(i,i^t):\ \overline{T}_{\Phi}(K) \leftarrow B\rightarrow \overline{T}^t_{\Phi^t}(K)\,.
\eqno(0.3)
$$
In this context, it turns out that
the fibers of the dual mirror fibrations lie in the monodromy invariant
homology class $\tau$ of middle dimension, as in the strict Calabi--Yau case.

\medskip

{\bf 0.3. Commutative and non--commutative tori and theta functions.}
In the second section of this paper we address the problem
of constructing ``quantum abelian varieties''. As we already mentioned,
what we actually construct is a linear space of
quantized theta functions (on a noncommutative torus 
with a period lattice)
endowed with partial multiplication. We regard this semiring 
as a quantum deformation of the universal multigraded
ring $\oplus \Gamma (L),\, L\in\roman{Pic}\,\Cal{A},\,[-1]^*(L)\cong L$,
which makes sense for any abelian variety $\Cal{A}$ (or indeed for
any projective variety).

\smallskip

Since the whole construction is algebraic, it
can be performed over any complete normed field, for example,
$p$--adic field, and applied to the classical
abelian varieties as well. Only those $p$--adic
abelian varieties admit a multiplicative
uniformization which have maximally degenerate stable reduction
modulo $p$. This is the definition of
the $p$--adic cusp. Diagrams of the type (0.2)
make full sense in this context as well, but $i^t$
admits no straightforward interpretation as anything
like symplectic form. It would be interesting
to investigate the meaning of $p$--adic $B$--fields
for strict Calabi--Yau manifolds.

\smallskip

We will now briefly discuss issues of non--commutative
geometry involved in our construction of quantized
theta functions. The standard approach is via
deformation of classical function rings, and this intuition
guided our initial construction in [Ma1]
and to a certain degree its extension presented in
this paper.
 
\smallskip

A complementary paradigm, made explicit
on many occasions in Connes' papers and the book [Co]
is the natural appearance of non--commutative rings
as objects encoding commutative spaces ``with bad geometric
properties'', typically quotients of commutative spaces
by non--separated equivalence relations.

\smallskip

For example, in the context of multiplicative uniformization
of abelian varieties $\Cal{A}\,(\bold{C})=T(\bold{C})/B$,
the multiplicative period lattice $B$
can degenerate also by ceasing
to be discrete, although keeping its rank constant,
and in this case it becomes natural to interpret $T(\bold{C})/i(B)$ in the realm of non--commutative geometry.

\smallskip

To be concrete, consider the case of one--dimensional
commutative torus $T.$ Its maximal toric compactification is $\bold{P}^1=\bold{G}_m\cup
\{0\}\cup \{\infty\}$. Choose $B=\bold{Z}$
so that $\Cal{A}=\bold{C}^*/(q^{\bold{Z}})$. 
On the level of diagram (0.3) we can choose here any $q\in\bold{P}^1(\bold{C})$.
When $|q|\ne 0, 1, \infty$, we get an elliptic curve $\Cal{A}=E_q$
fibered by images of real tori $|z|=\roman{const}$ (where $z$
is the coordinate on $\bold{C}^*$) over the circle $\bold{R}/\bold{Z}
\,\roman{log}\,|q|$. Dualizing this fibration, we will get the mirror
dual elliptic curve. 
The points $0, \infty$ become divisors on a modular
curve and provide the familiar
degeneration picture in algebraic geometry. Points where $q$
is a root of unity also become visible in algebraic geometry
as cusp points of modular curves of higher levels.
However, $q$ of infinite order with $|q|=1$ are not considered in algebraic
geometry at all. One way to interprete the space
$\bold{G}_m/(q^{\bold{Z}})$ in this case is to identify it with
(one of the versions of) the {\it two--dimensional} non--commutative
torus $T_q$, whose function ring is generated by $x^{\pm 1},y^{\pm 1}$ satisfying the
commutation relation $xy=qyx.$ See Appendix  for an informal
explanation of this in the context of noncommutative geometry
\`a la Alain Connes. 

\smallskip

One remarkable trace of this
origin of $T_q$ for $q=e^{2\pi i\tau},\,\tau\in \bold{R},$ as a limiting elliptic curve is the re--appearance
of the modular group $SL(2,Z)$ as a symmetry group
in the non--commutative situation: acting on $\tau$, it produces
Morita equivalent function rings of $T_q$'s. See [RiS] for a thorough
discussion of this in arbitrary dimension.

\smallskip

What I want to stress here, is a somewhat neglected
complementary aspect of this picture: namely, that even
if the group $q^{\bold{Z}}$ is discrete, 
$T_q$ still can be viewed as a legitimate
incarnation of the elliptic curve $E_q$ in the non--commutative world.
A systematic treatment of the correspondence
between, say, coherent sheaves on $E_q$ and modules over $T_q$
remains a problem for future, but see the recent preprint [BEG] for
some precise facts about this correspondence, and
many interesting suggestions are contained in [So].

\smallskip

This remark throws some light on the problems left open
with our approach to quantized theta--functions,
for example, that of functional equations
corresponding to a change of cusp. Some of our
non--commutative abelian varieties $T(H,\alpha )/B$
where $T(H,\alpha )$ is a {\it non--commutative torus}
with quantization parameter $\alpha$ (see 2.1)
can be understood as a result of taking a quotient
of a {\it commutative torus } $T(H^{\prime}, 1)$ of
halved rank, by a {non--discrete period subgroup with too many
generators}.
This might provide a bridge between our
construction and that of Weinstein ([We]) remaining
entirely in the realm of commutative geometry.  

\medskip

{\it Acknowledgement.} I am grateful to D.~Orlov
for illuminating correspondence on Abelian mirrors
and comments on the paper [GolLO].
I learned from Y.~Soibelman the philosophy
of treating boundaries of various moduli spaces
as bridges to the non--commutative realm.
A.~Polishchuk has drawn  my attention to the paper [Fu]
and sent me a copy of it.

\bigskip

\centerline{\bf 1. Mirror symmetry for complex tori and abelian varieties}

\medskip

{\bf 1.1. Toric formalism.} Since we will have to work
with several different types of tori which must
be carefully distinguished, we start
with some terminological conventions. 

\smallskip

Let $K$ be a ground field,
$H$ a free abelian group of rank $n$ which we will always write
additively. {\it Algebraic torus $T$ with character group $H$}
is the affine  spectrum of the group ring of 
$H$, and $H=\roman{Hom}\, (T, \bold{G}_m)$.
Put $H^t:=\roman{Hom}\,(\bold{G}_m,T)$. Groups $H, H^t$ are
connected by the canonical duality map $H\times H^t\to
\roman{Aut}\,\bold{G}_m=\bold{Z}.$ 
\smallskip

An element $h\in H$ considered as a character of $T$ will be denoted
$e(h)$ ($e$ for exponential). We have $e(h+h^{\prime})=e(h)e(h^{\prime}).$ 

\smallskip

In the next section, we will consider noncommutative algebraic tori
as well, for which the multiplication rule for $e(h)$
is  twisted: see (2.3) below.

\medskip

{\bf 1.1.1. Duality of algebraic tori.} If $T$ is an algebraic torus
with character group $H_T$, the dual algebraic torus $T^t$
has the character group $H_{T^t}:=H_T^t =\roman{Hom}\,(H_T,\bold{Z})$.

\medskip

{\bf 1.1.2. Periods and abstract tori.} Let $H$ be the
character group of $T$. Denote by
$B$  another free abelian group of the same rank $n$, and by
$i:\,B\to T(K)$ a homomorphism which we will call ``period map''.
Since $T(K)=\roman{Hom}\,(H,K^*)$, to give $i$ is the same as to give a pairing
$H\times B\to K^*$ such that $(h,b)$ is the value of
$e(h)$ at the point $i(b)$.

\smallskip

We will refer to the quotient space $T/i(B)$ as 
{\it an abstract torus}, and to $T$ as {\it its covering algebraic
torus}. One may imagine $T/i(B)$ simply
as a functor of points on $K$--algebras $R\mapsto T(R)/i(B).$ 
In the case when $K=\bold{C}$ and $i$ is an injection with discrete image,
$T(\bold{C})/i(B)$ is {\it a complex torus}, which may admit
a structure of abelian variety (it is then unique).
However, a large part of our elementary formalism will
not depend on additional assumptions. Thus, as explained in the Introduction,
we will be able to include into our mirror picture
$p$--adic tori and abelian varieties, and eventually
non--commutative tori viewed as models of $T(\bold{C})/i(B)$
in non--commutative geometry.
\smallskip

In order not preclude the eventual interpretation
of the space $\Cal{A}=T/i(B)$ we will consistently identify
$\Cal{A}$ with the triple 
$$
(H_{\Cal{A}},B_{\Cal{A}}, (\,,)_{\Cal{A}}:\,
H_{\Cal{A}}\times B_{\Cal{A}}\to K^*)
\eqno(1.1)
$$ 
as above. 

\medskip

{\bf 1.1.3. Poincar\'e dual abstract tori.} By definition,
Poincar\'e duality interchanges characters and periods.
More precisely, $\Cal{A}$ and $\widehat{\Cal{A}}$ are Poincar\'e
dual if
$$
H_{\widehat{\Cal{A}}}=B_{\Cal{A}},\
B_{\widehat{\Cal{A}}}=H_{\Cal{A}},\
(b,h)_{\widehat{\Cal{A}}}= (h,b)^{-1}_{\Cal{A}}\ \roman{for}\
h\in H_{\Cal{A}},\,b\in B_{\Cal{A}}.
\eqno(1.2)
$$
For abelian varieties this agrees
with the classical definition. Notice however, that
a choice of the covering algebraic torus is an additional structure, and we
explicitly extend Poincar\'e duality to this context.

\medskip

{\bf 1.1.4. Framed tori.} {\it A framing} of the abstract 
torus (1.1) is a map $i^t:\,B_{\Cal{A}}\to T_{\Cal{A}}^t(K)$.
{\it A framed abstract torus} is a pair $(\Cal{A},i^t)$.

\smallskip

A framing of the complex torus or abelian variety consists
of its representation as an abstract torus and framing
of that abstract torus.

\medskip

{\bf 1.1.5. Mirror dual framed abstract tori.} Two
framed abstract tori $(\Cal{A},i^t)$ and $(\Cal{B},i)$ 
are called {\it mirror dual}, or {\it mirror partners}, if
their covering algebraic tori are dual, and their periods are
explicitly identified. More precisely, the relation of
mirror partnership is provided by diagrams of the form
(0.2). If one thinks about $i,i^t$ in terms of the respective
character/period pairings (1.1), the mirror duality
is provided by pairings
$$
H_{\Cal{A}}\times B_{\Cal{A}}\to K^*,\
H^t_{\Cal{A}}\times B_{\Cal{A}}\to K^* .
\eqno(1.3)
$$
A framing is called {\it non--degenerate} if both kernels of
the respective pairing are trivial.
\smallskip

This notion of mirror duality is the main definition of this section. We start
with studying it for $K=\bold{C}.$

\medskip

{\bf 1.2. Complex tori.} Assume
that $i(B)$ is discrete in $T(\bold{C})$. Put
$$
\Gamma =\Gamma_{\bold{Z}}=\pi_1(T(\bold{C})/i(B),0)=H_1(\Cal{A},\bold{Z}),
$$ 
$\Gamma_{\bold{R}}=\Gamma_{\bold{Z}}\otimes\bold{R}$
and similarly for $\Gamma_{\bold{C}}$. If $T(\bold{C})/i(B)$
is denoted $\Cal{A}$, we may write $\Gamma_{\Cal{A}},$
$\Gamma_{\Cal{A},\bold{R}}$,
etc. 

\smallskip

The real space  $\Gamma_{\bold{R}}$ can be identified
with the Lie algebra of $\Cal{A}$, and the exponential
map $\roman{exp}:\, \Gamma_{\bold{R}}\to \Cal{A}$
with kernel $\Gamma_{\bold{Z}}$ is the universal covering
of $\Cal{A}$. 

\smallskip

Let $H=H_{\Cal{A}}$ be the character group of $T$. The map
$h\mapsto \dfrac{1}{2\pi i}\,\dfrac{d\,e(h)}{e(h)}$
induces canonical identification $H=H^1(T(\bold{C}),\bold{Z})$.
Hence we have an exact sequence
$$
0\to H_{\Cal{A}}^t\to H_1(\Cal{A},\bold{Z})\to B_{\Cal{A}}\to 0 
\eqno(1.4)
$$
where the third arrow is induced by $i$. Similarly, we
have for a mirror dual framed torus $(\Cal{B}=T^t(\bold{C})/i^t(B),i)$
$$
0\to H_{\Cal{A}}\to H_1(\Cal{B},\bold{Z})\to B_{\Cal{A}}\to 0 
\eqno(1.5)
$$
where the third arrow is now induced by $i^t$.

\medskip

{\bf 1.2.1. Mirror partners as dual real torus fibrations.}
In this subsection we will show that our definition
of mirror partners $(\Cal{A},\Cal{B})$  
over $\bold{C}$ naturally fits into the general context
of Lagrangian/complex duality: see [AP], [StYZ], [Gr2].

\smallskip

We start with a brief description of this context.

\smallskip

Let  $(X,\omega )$ be a $C^{\infty}$ symplectic
manifold, endowed with a submersion
$p_X:\,X\to U$ whose fibers are Lagrangian tori.
We will fix also a Lagrangian section $0_X:\,U\to X.$

\smallskip

Using $\omega$, we can identify the bundle of Lie algebras
of the tori $p^{-1}_X(u),\,u\in U,$ with the cotangent bundle
$T^*_U.$ Hence we have a canonical isomorphism
$X=T^*_U/H$ where $H$ is a Lagrangian sublattice in $T^*_U$
with respect to the lift of $\omega$ which is the standard
symplectic form on the cotangent bundle. There exists
also a canonical flat symmetric connection
on $T^*_U$ for which $H$ is horizontal.

\smallskip

The local system $H^t=\Cal{H}om\,(H,\bold{Z})$  is embedded as
a sublattice into $T_U$, and we can define {\it the mirror
partner} of $(p_X:\,X\to U,\omega ,0_X)$ as the toric fibration
$Y:=T_U/H^t$ endowed with the projection  to the same
base $p_Y:\,Y\to U$
and the zero section $0_Y$.

\smallskip

Passing from
$X$ to $Y$ we have lost the symplectic form. To compensate for this loss,
we have acquired a complex structure $J:\,T_Y\to T_Y$
which can be produced from $(p:\,X\to U,\omega ,0_X)$
in the following way. The flat connection on
$T_U$ obtained by the dualization from $T^*_U$
produces a natural splitting $T_Y=p_Y^*(T_U)\oplus p_Y^*(T_U).$
With respect to this splitting, $J$ acts as $(t_1,t_2)\mapsto (-t_2,t_1).$

\smallskip

Conversely, suppose that we have a complex manifold $Y$ endowed
with a fibration by real tori $Y\to U$ with zero section,
such that the operator of complex structure along the zero section
identifies $T_U$ with the bundle of Lie algebras of fibers.
Then we can consecutively construct the lattice $H^t\subset T_U$,
the dual fibration $X:=T^*_U/H$ and the symplectic form
on $X$ coming from the cotangent bundle.

\smallskip

Now we can return to complex tori.

\smallskip

Put $S^1=\{ |z|=1\,|\,z\in\bold{C}\}.$ We have
the Lie group isomorphism $\bold{C}^*\to S^1\times\bold{R}:$
$z\mapsto (z/|z|, \roman{log}\,|z|)$. This induces
an isomorphism
$$
(\alpha ,\lambda ):\ T(\bold{C})\to \roman{Hom}\,(H,S^1)
\times \roman{Hom}\,(H, \bold{R}).
\eqno(1.6)
$$
If $i(B)$ is discrete of maximal rank which I will assume,
then $\lambda\circ i(B)$ is an additive lattice in the real space
$\roman{Hom}\,(H, \bold{R})$.
Thus (1.6) produces a real torus fibration of $T(\bold{C})$ 
over the base which
is as well a real torus of the same dimension:
$$
0\to \roman{Hom}\,(H,S^1)\to T(\bold{C})/i(B)\to
\roman{Hom}\,(H,\bold{R})/\lambda\circ i(B)\to 0 \,.
\eqno(1.7)
$$
Similarly, we have 
$$
0\to \roman{Hom}\,(H^t,S^1)\to T^t(\bold{C})/i^t(B)\to
\roman{Hom}\,(H^t,\bold{R})/\lambda^t\circ i^t(B)\to 0 
\eqno(1.8)
$$
where $\lambda^t$ is defined for $T^t$ in the same way
as $\lambda$ for $T$. Let us identify linear real spaces $H_{\bold{R}}$
with $H^t_{\bold{R}}$ in such a way that lattice points
$\lambda\circ i(b)$ and $\lambda^t\circ i^t(b)$ are identified
for all $b\in B$. Then (1.7) and (1.8) become
dual real torus fibrations over the common base.

\smallskip

The relevant complex structures in our context come from covering tori.
We have to introduce  symplectic forms.
Let us construct, say, $\omega_{\Cal{A}}$. From (1.8)
one sees that $\Cal{A}=T(\bold{C})/i(B)$ can be obtained
as quotient space of the tangent bundle of the base
by a lattice. The tangent bundle (and the lattice)
is canonically trivialized, and its fiber is $H_{\bold{R}}.$
Using the two framings, we have identified $H_{\bold{R}}$ 
with $H^t_{\bold{R}},$ that is, tangent bundle with cotangent bundle.
The canonical symplectic form on the cotangent bundle becomes our
$\omega_{\Cal{A}}.$ Clearly, fibers of (1.7)
are Lagrangian tori. It remains to check that $\omega_{\Cal{A}}$
determines $I_{\Cal{B}}$ as above, but this is quite straightforward.

\medskip

{\bf 1.2.2. Maximal degeneration point and monodromy.}
Let us consider now the situation, described
by (one half of) the diagram (0.3). More precisely,
consider the space of maps of $B$ to a neighbourhood
of a point of maximal degeneration in some
toric compactification $\overline{T}_{\Phi}.$
Such a point is a zero--dimensional orbit of $T$,
thus it corresponds to a maximal cone in $\Phi$.
Assume for simplicity that it is the simplicial cone generated
by a basis of $H^t$. This means that we identify
$H^t$ with $\bold{Z}^n$, $T(\bold{C})$ with $(\bold{C}^*)^n$,
and choose as partial compactification the imbedding
$(\bold{C}^*)^n\subset \bold{C}^n$. Let $D_r$
be the $r$--th coordinate hyperplane in $\bold{C}^n$.

\smallskip

Choosing a basis of $B$ as well, that is, identifying
it with $\bold{Z}^n$, we see that the region of the
partially compactified moduli space of 
multiplicatively uniformized complex tori that we are interested in can be
identified with an open subspace of the matrix space
$\bold{C}^{n\times n}$ whose columns generate
a multiplicative sublattice. The discriminant
locus in this region consists of its intersection
with $\cup D_{rs}$ where $D_{rs}$ is 
the $r$--th coordinate hyperplane in the $s$--th copy
of $\bold{C}^n$ times other copies. The origin
(intersection of all $D_{rs}$) is the maximum degeneration
point.

\smallskip

Let $q=(q_{rs})$ be a point of the moduli space, $\Cal{A}_q$
the respective torus. Denote by $\gamma_r\in H_1(\Cal{A}_q,\bold{Z})$
the image of the $r$--th $S^1$ in (1.7) with counterclockwise
orientation. Let $M_{rs}$ be the monodromy action
of a small counterclockwise loop around $D_{rs}$ in the
moduli space. All cycles $\gamma_r$ are monodromy invariant.
Let $\beta_s$ be any lift to $H_1(\Cal{A}_q,\bold{Z})$
of the $s$--th 1--cycle in the base torus in (1.7).
Then $M_{rs}$ transforms $\beta_s$ into $\beta_s+\gamma_r$
and leaves other 1--cycles invariant.

\smallskip

Thus, the homology class $\tau$ of any fiber of (1.7)
generates the cyclic group of invariant cycles
of middle dimension. 

\smallskip

This statement holds independently of the choice of
the simplicial fan and a maximal degenerating cone in it.
In this sense, the choice of a covering torus alone encodes essential
information about large complex structure.

\medskip

{\bf 1.3. Framings and well--becoming pairs.} In this subsection
we will compare our construction of mirror partners
with that of [GolLO]. In that paper, the additional structure
on $\Cal{A}$ is a complex--valued 2--form $\omega$
rather than symplectic form as above. I will show that
this as well can be  related to an appropriate framing. 

\smallskip

Consider first 
a pair of framed complex tori $(\Cal{A},i^t)$
and $(\Cal{B},i)$ which are mirror partners as above.
In [GolLO], sec. 10,
the authors use  decompositions 
$$
\Gamma_{\Cal{A}}=\Gamma_{1,\Cal{A}}\oplus\Gamma_{2,\Cal{A}},\
\Gamma_{\Cal{B}}=\Gamma_{1,\Cal{B}}\oplus\Gamma_{2,\Cal{B}} .
$$
We will call such  decompositions {\it compatible} with
our choice of covering tori, if $\Gamma_{1,\Cal{A}}=H^t_{\Cal{A}},$
$\Gamma_{1,\Cal{B}}=H_{\Cal{A}}$ as in (1.4), (1.5).
Thus compatible decompositions are simply splittings
of (1.4) and (1.5). The spaces $\Gamma_{\Cal{A},\bold{R}},$
$\Gamma_{\Cal{B},\bold{R}}$
are endowed respectively with 
complex structures $I_{\Cal{A}},I_{\Cal{B}}$ coming from
covering tori. Clearly, 
$$
\Gamma_{\Cal{A},\bold{R}} =H^t_{\Cal{A},\bold{R}}\oplus
I_{\Cal{A}}H^t_{\Cal{A},\bold{R}}
$$ 
so that compatible splittings satisfy
the condition 10.3.1 (2) of [GolLO]. 

\smallskip

Consider now only $\Cal{A}$, but equipped with a class
$\omega\in H^2(\Cal{A},\bold{C})$ interpreted as an antisymmetric
complex--valued form on $\Gamma_{\Cal{A}}$. Assume that the extension of $\omega$
to $\Gamma_{\Cal{A},\bold{R}}$ is $I_{\Cal{A}}$--invariant (see [GolLO], last lines
of 1.4 for an explanation of this condition).
Assume moreover that there exists a compatible
decomposition of $\Gamma_{\Cal{A}}$  such that $\Gamma_{1,\Cal{A}}$ and $\Gamma_{2,\Cal{A}}$
are $\omega$--isotropic. This means that $(\Cal{A},\omega )$ is
a well--becoming pair in the sense of [GolLO], 10.3.1.
Looking at (1.4), we see that $\omega$ can be uniquely reconstructed from
its restriction which we also denote $\omega$ 
$$
\omega:\, H^t_{\Cal{A}}\times B_{\Cal{A}}\to \bold{C} .
\eqno(1.9)
$$
Exponentiating (1.9) we produce a framing $i^t:\,B\to T^t(\bold{C})$ in the form
(1.3), that is
$$
e(h)(i^t(b))=e^{2\pi i\,\omega (h,b+H^t_{\Cal{A}})},\ h\in H^t_{\Cal{A}},\,b\in B.
\eqno(1.10)
$$ 
Clearly, $i^t$ remains the same, if we choose another compatible isotropic splitting. Moreover, it does not change if we add to $\omega$
another pairing taking values in $2\pi i\bold{Z}$.

\smallskip

In this way we get a map from the set of all well--becoming
pairs $(\Cal{A},\omega )$ admitting compatible isotropic splittings
in the sense of [GolLO] 
to the set of framed abstract tori $(\Cal{A},i^t)$ with non--degenerate
framings in our sense.

\smallskip

We can now complete the comparison of our definition
of mirror duality with that of [GolLO]. 

\medskip

\proclaim{\quad 1.3.1. Theorem} Let $(\Cal{A},i^t)$,
$(\Cal{B},i)$ be a mirror dual pair of framed complex abstract tori,
admitting lifts $(\Cal{A},\omega_{\Cal{A}})$, $(\Cal{B},\omega_{\Cal{B}})$
to well--becoming pairs.
Then $(\Cal{A},{\omega}_{\Cal{A}})$, $(\Cal{B},{\omega}_{\Cal{B}})$ 
are mirror dual in the sense
of [GolLO].
\endproclaim

\smallskip

{\bf Proof.} We will compare our setting with that of [GolLO],
10.4 and 10.4.1. Groups $\Gamma_1,\Gamma_2$ in [GolLO]
are our $H^t_{\Cal{A}}, B$, fundamental group of the
mirror dual torus is $\Gamma^t_1\oplus\Gamma_2$. This means
that the {\it real} covering torus of their $\Cal{B}$ is the same
as ours, that is  $T^t$. It remains to compare the complex structures.
In our case it is simply induced from $T^t(\bold{C}).$
In [GolLO] it is described with the help of the complex
structure operator $I$
acting upon  $(\Gamma^t_1\oplus\Gamma_2)\otimes\bold{R}$ 
produced from $\omega$
in two steps: via formula (14) and subsequent projection
described in 10.4. 

\smallskip

In order to check that they coincide, we will reproduce a part
of the argument in [GolLO] in the form which hopefully clarifies
the meaning of their crucial formula (14).

\smallskip

Since the following construction must be considered
in two different situations, described in [GolLO] 9.2
and 10.4 respectively, we slightly change the scope
of our notation. From now on, $\Gamma_1, \Gamma_2$
will denote two abstract free abelian groups of the same
finite rank, $\Gamma = \Gamma_1\oplus \Gamma_2$,
$\Gamma^{\prime} = \Gamma_1^t\oplus \Gamma_2$.
 Introduce
two real tori
$$
\Cal{C}=(\Gamma_1\oplus \Gamma_2)_{\bold{R}}/\Gamma_1\oplus \Gamma_2,\
\Cal{C}^{\prime}=(\Gamma_1^t\oplus \Gamma_2)_{\bold{R}}/\Gamma_1^t\oplus \Gamma_2 .
\eqno(1.11)
$$
Consider the data of two types.

\smallskip

(i) Complex structures on $\Cal{C}^{\prime}$ described
by the operators $I$ on $(\Gamma_1^t\oplus \Gamma_2)_{\bold{R}}$
such that their crossover components 
$$
I_{12}:\, \Gamma_{2,\bold{R}}\to\Gamma_{1,\bold{R}}^t,\
I_{21}:\, \Gamma_{1,\bold{R}}^t\to\Gamma_{2,\bold{R}}
$$
are bijective.

\smallskip

(ii) Forms $\omega\in \wedge^2\Gamma^t_{\bold{C}}$ for which
$\Gamma_{1,\bold{C}}$ and $\Gamma_{2,\bold{C}}$ are maximal
isotropic.

\smallskip

We will establish a bijection between them in the following way.

\smallskip

Let us start with a complex structure $I$ in $\Cal{C}^{\prime}$.
It determines (and is determined by) the space of invariant
holomorphic 1--forms on $\Cal{C}^{\prime}$. Integrating
them over $\Gamma_1^t\subset H_1(\Cal{C}^{\prime},\bold{Z})$,
we will get all additive maps $\Gamma_1^t\to \bold{C}$,
in particular, all elements of $\Gamma_1$. So we have an embedding
$\Gamma_1\to H^0(\Cal{C}^{\prime},\Omega^1):\ \gamma\to\nu_{\gamma}$
such that for all $\beta\in\Gamma_1^t, \gamma\in\Gamma_1,$
$$
(\beta ,\gamma )=\int_{\beta}\nu_{\gamma} .
\eqno(1.12)
$$
This allows one to define a non--degenerate scalar product
$\langle\, ,\rangle :\Gamma_2\otimes\Gamma_1\to\bold{C}$:
$$
\langle\gamma_2 ,\gamma_1 \rangle =\int_{\gamma_2}\nu_{\gamma_1} .
\eqno(1.13)
$$
Finally, we can extend it to a complex skew--symmetric form $\omega$
on $\Gamma$ declaring $\Gamma_1$ and $\Gamma_2$ to be isotropic:
$$
\omega ((\gamma_1,\gamma_2), (\gamma_1^{\prime},\gamma_2^{\prime}))
=\langle\gamma_2 ,\gamma_1^{\prime} \rangle
- \langle\gamma_2^{\prime} ,\gamma_1 \rangle .
\eqno(1.14)
$$
If we choose a basis of $\Gamma_1^t,\Gamma_2$ and a basis
of holomorphic 1--forms whose period matrix over $\Gamma_1^t$
is the identity $E$, then the Gram matrix of the pairing $\langle\, ,\rangle$
will be just the second half of the total period matrix. Let us
denote it $\tau$. Then $e^{2\pi i\,\tau}$ is the matrix
generating the multiplicative period lattice in the covering
complex torus $\Gamma_{\bold{C}}^{\prime}/\Gamma_1^t$ which in view
of (1.14) agrees with (1.10).

\smallskip

Arguing now in reverse direction, we will show that knowing
$\tau$ we can reconstruct the operator $I$ in the same basis
and get essentially the [GolLO] formula.
In fact, $I$ is uniquely determined by the requirement
that for all $\gamma\in \Gamma^{\prime}$ and all holomorphic $\nu$
we have $\int_{I\gamma}\nu =i\int_{\gamma}\nu .$ Hence
to find $I$ we must solve the matrix system
$$
(\roman{Re}\,\tau +i\,\roman{Im}\,\tau ,\,E)\,
\pmatrix X& Y\\U& V\endpmatrix =
(-\roman{Im}\,\tau\,+ i\roman{Re}\,\tau ,\,iE)
$$
which gives
$$
I=\pmatrix X& Y\\U& V\endpmatrix
=\pmatrix 
(\roman{Im}\,\tau )^{-1}\, \roman{Re}\,\tau &
(\roman{Im}\,\tau )^{-1}\\
-\roman{Im}\,\tau - \roman{Re}\,\tau\, 
(\roman{Im}\,\tau )^{-1}\,\roman{Re}\,\tau &
-\roman{Re}\,\tau 
\,(\roman{Im}\,\tau )^{-1}\endpmatrix
\eqno(1.15)
$$
For the first application of (1.15), let us choose a real torus $\Cal{A}$ with period lattice
$\Gamma_{\Cal{A}}$ and put $\Cal{C}= \Cal{A} \oplus\widehat{\Cal{A}}$
where $\widehat{\Cal{A}}$ is the Poincar\'e dual real torus, that is
$\Gamma_{\widehat{\Cal{A}}}=\Gamma^t_{\Cal{A}}$ (this does not contradict
our multiplicative description (1.2) although it may not be
immediately obvious). Put $\Gamma_1=\Gamma_{\Cal{A}},$ $\Gamma_2=
\Gamma^t_{\Cal{A}}.$ Comparing our formula (1.15) with [GolLO] (14)
in this situation, we see that our $I$ coincides with their
$I_{{\omega}}$. In [GolLO], 8.4, $\Cal{A}$ additionally possesses a complex
structure, which produces the canonical complex structures
on $\widehat{\Cal{A}}$, $\Cal{C}$ and $\Cal{C}^{\prime}$, say, $J.$
Moreover,
$\omega$ is restricted to lie in the complexified
N\'eron--Severi group, and as a result $I$ commutes with the
inherited complex structure $J$ on $\Cal{C}^{\prime}.$

\smallskip

The setup which we are discussing in the Theorem 1.3.1 is
that of [GolLO] 10.4 and 10.5. The relevant torus $\Cal{C}$
is now $\Cal{A}$, and its homology lattice is 
now split by the choice of a compatible decomposition like
in 1.3. Extension of this splitting to $\Cal{A} \oplus\widehat{\Cal{A}}$
produces the period matrix $\tau$ which is block diagonal
and consists of two blocks. Formula (1.15) is still valid
for the mirror complex structure, when one replaces $\tau$
in it by the respective block. Putting everything
together we see that our complex structure determined
essentially by (1.10) indeed agrees with that of [GolLO].  
\medskip

{\bf 1.3.2. Remark.} [GolLO] contains several tentative descriptions
of mirror dual pairs differing mostly by the exact choice of
the additional structure that should be added to $\Cal{A}$.
Our Theorem 1.3.1 together with the
Theorem 10.5 in [GolLO] indicates that the notion of a 
well--becoming pair endowed with a choice of one half
of an isotropic decomposition captures just right amount of information.
Replacing this structure by that of framing, we  make explicit the important
aspect of ``large complex structure'' in the case $K=\bold{C}$  
and simultaneously extend
the definition to abstract tori over arbitrary fields.

\bigskip

\centerline{\bf 2. Quantized theta--functions and abelian varieties.}

\bigskip

{\bf 2.1. Category of non--commutative tori.} 
Let $H$ be a free abelian group of finite rank and
$\alpha :\,H\times H\to K^*$ an alternating pairing: for all $h,g\in H$
$$
\alpha (h,g)=\alpha (g,h)^{-1},\quad \alpha (h_1+h_2,g)=
\alpha (h_1,g) \alpha (h_2,g).
\eqno(2.1)
$$
A morphism $f:\,(H_1,\alpha_1)\to (H_2,\alpha_2)$ is a group homomorphism
$f:\,H_1\to H_2$ such that for all $h,g\in H_1$ we have
$$
\alpha_2^2(f(h),f(g))=\alpha_1^2(h,g).
\eqno(2.2)
$$
The bilinear form 
$$
\varepsilon_f (h,g):= \alpha_1(h,g)\alpha_2^{-1}(f(h),f(g))
$$
with values in $\{\pm 1\}$ is called {\it the characteristic of $f$}
(and of $F$).

\smallskip

Any such pair $(H,\alpha )$ will be called {\it the character group}
of the  non--commutative torus $T(H,\alpha )$ whose
{\it ring of algebraic functions} $Al(H,\alpha )$ is the linear space spanned
over $K$ by the symbols $e(h),\,h\in H,$ with multiplication rule
$$
e(h) e(h^{\prime}) =\alpha (h,h^{\prime})\,e(h+h^{\prime}).
\eqno(2.3)
$$
We may write $e_{H,\alpha}(h)$ for $e(h)$ if need be.

\smallskip
Notice that $\varepsilon (h):=\alpha (h,h)$ is a character
of $H$ taking values $\pm 1$, and that from (2.3) we get the following
formulas:
$$
e(h_1) e(h_2) e(h_3)= \alpha (h_1,h_2)\,\alpha (h_1,h_3)\,
\alpha (h_2,h_3)\, e(h_1+h_2+h_3),
\eqno(2.4)
$$
$$
e(h)^{-1}=\varepsilon (h)\,e(-h).
\eqno(2.5)
$$
We can also consider the two-sided $Al\,(H,\alpha )$--module
of {\it formal functions} $Af\, (H,\alpha )$ consisting
of infinite linear combinations $\sum_h a_he(h),\,a_h\in K$,
and, in the case of a complete normed field $K$ and {\it an unitary} 
quantization parameter $\alpha$ (that is, $|\alpha |=1$) the ring
of {\it analytic functions} $An\, (H,\alpha )$ consisting
of those formal functions for which $|a_h|\,\|h\|^N\to 0$
for any $N$ as $\|h\|\to \infty$, $\|h\|$ being any
Euclidean norm on $H$.

\smallskip

The form $\alpha$ can be called {\it the quantization parameter.}
When $\alpha \equiv 1$, we get the usual notions of commutative geometry,
so that $T(H,1)$ is the algebraic torus with character group $H$.

\smallskip

{\it A morphism} $F:\,T(H,\alpha )\to T(H^{\prime},\alpha^{\prime})$,
by definition, is given by the contravariant $K$--algebra homomorphism
$F^*:\,Al\,(H^{\prime},\alpha^{\prime})\to Al\,(H,\alpha )$.

\medskip

\proclaim{\quad 2.1.1. Proposition} a) The set of invertible elements
of $Al\,(H,\alpha )$ is $\{ a\,e(h)\,|\,a\in K^*, h\in H\}.$
If  $F:\,T(H_2,\alpha_2 )\to T(H_1,\alpha_1)$ is a morphism of non--commutative tori, then the induced map $f=[F]:\, H_1\to H_2$
determined by $F^*(e(h))=a_he(f(h)),\, a_h\in K^*,$ 
is additive and satisfies (2.2)
and thus is a morphism of character groups.

\smallskip

b) The set of all morphisms $F$ with fixed $f=[F]$ is either empty, or
has a natural structure of the principal homogeneous space
over the group $T(H_1,1)(K)=\roman{Hom}\,(H_1,K^*).$
In particular, if the characteristic of $f$ is 1,
then $F^*:\,e(h)\mapsto e(f(h))$ extends to a uniquely
defined morphism of rings of algebraic functions.

\smallskip

c) Any morphism $F^*$ extends to $Af$ by $F^*(\sum a_h e(h))=
\sum a_hF^*(e(h)).$ If $K$ is normed and $\alpha$ unitary,
then this extension maps analytic functions to analytic.
\endproclaim

\smallskip

{\bf Proof.} The first statement follows from the fact that $H$ can 
be endowed with the structure of a well--ordered group.
For any such structure, the highest (resp. lowest) terms
of a product are products of the highest (resp. lowest) terms,
so that an invertible element coincides with its highest
and lowest term.

\smallskip

To prove the second statement, rewrite the equality
$$
F^*(e(h)\,e(g))=F^*(e(h))\,F^*(e(g))
$$
using (2.3).
Comparing the $e$--terms, we see that $f(h+g)=f(h)+f(g).$
Comparing the coefficients, we get
$$
a_h a_g a^{-1}_{h+g} = 
\alpha_1(h,g)\alpha_2^{-1}(f(h),f(g)).
\eqno(2.6)
$$
The left hand side is a symmetric form in $h,g$, whereas the
right hand side is alternate. Therefore
this form takes values $\{\pm 1\}$. Hence (2.2) holds,
and (2.6) is the characteristic $\varepsilon_f (h,g)$ of $f$.

\smallskip

Finally, let $f$ be a morphism of character groups with
characteristic $\varepsilon$. Then ring morphisms $F^*$
with $[F]=f$ bijectively correspond to the solutions
$\{a_h\,|\,h\in H_1\}$ of the equations 
$a_h a_g a^{-1}_{h+g} =\varepsilon_f (h,g)$. If one such
solution exists, then all others are of the form
$a_hc(h)$ where $c:\,H_1\to K^*$ is an arbitrary homomorphism.

\smallskip

The remaining statements are straightforward.

\medskip

{\bf 2.1.2. Direct product.} By definition,
the ring of algebraic functions of $T(H_1,\alpha_1)\times T(H_2,\alpha_2)$
is the tensor product of the respective rings. We can write
$$
T(H_1,\alpha_1)\times T(H_2,\alpha_2)=
T(H_1\oplus H_2,\alpha_1\oplus \alpha_2)
$$
by identifying
$$
e_{H_1,\alpha_1}(h_1)\otimes e_{H_2,\alpha_2}(h_2) =
e_{H_1\oplus H_2,\alpha_1\oplus \alpha_2}((h_1,h_2)).
$$
\smallskip

{\bf 2.1.3. Some standard morphisms.}
{\it (i) Shifts.} Any point $b\in T(H,1)(K)=\roman{Hom}\,(H,K^*)$
determines an automorphism $b^*$ of $T(H,\alpha )$:
$$
b^*(e(h)):=h(b)e(h),
\eqno(2.7)
$$
where from now on we denote by $h(b)$ the value of $e_{H,1}(h)$
at the point $b$. 

\smallskip

{\it (ii) Multiplication by $n$.} This is the morphism
$$
[n]:\, T(H,\alpha )\to T(H,\alpha^{n^2})
$$
defined by
$$
[n]^*(e_{H,\alpha^{n^2}}(h))=e_{H,\alpha}(nh).
\eqno(2.8)
$$
For $n=-1$ it is an endomorphism of $T(H,\alpha )$.
It is also an endomorphism, if $\alpha$ takes values
in roots of unity of degree $d$ and $n^2\equiv 1\,\roman{mod}\,d.$

\smallskip

The commutation rule with shifts is
$$
b^*\circ [n]^*=[n]^*\circ (nb)^* .
\eqno(2.9)
$$
\smallskip

{\it (iii) External multiplication.} It is the morphism
$$
m_{\alpha ,\beta}:\, T(H,\alpha )\times T(H,\beta )\to
T(H,\alpha\beta )
$$
defined by
$$
m^*_{\alpha ,\beta} (e_{H,\alpha\beta}(h))=e_{H,\alpha}(h)\otimes 
e_{H,\beta}(h) .
\eqno(2.10)
$$

\smallskip

{\it (iv) Mumford's morphism.} This is the morphism
$$
M:\,T(H\oplus H,\alpha\oplus\alpha )\to T(H\oplus H,\alpha^2\oplus\alpha^2),
$$
$$
M^*(e(h,g))=e(h+g,h-g).
\eqno(2.11)
$$
It is well defined, because
$$
(\alpha\oplus\alpha )[(h+g,h-g),(h^{\prime}+g^{\prime}, h^{\prime}-g^{\prime})]
$$
$$
=\alpha (h+g, h^{\prime}+g^{\prime}) \alpha (h-g, h^{\prime}-g^{\prime})
=\alpha^2(h,h^{\prime}) \alpha^2(g,g^{\prime}).
$$
\medskip

{\bf 2.2. Periods.} We choose and fix an abelian
group  of periods $B\subset T(H,1)(K).$ 
The period group is
written additively; it acts upon $T(H,\alpha )$ by shifts.

\smallskip

Trying to make sense of the quotient
$T(H,\alpha )/B$ we will study formal or analytic
functions on $T(H,\alpha )$ with automorphic properties
with respect to the the group $\{b^*\,|\,b\in B\}.$

\medskip

\proclaim{\quad 2.3. Definition} A (two--sided) theta multiplier $\Cal{L}$
for the non--commutative torus $T(H,\alpha )$ and period group $B$
consists of the data $\Cal{L} = (h_l,h_r,\psi , (\,,))$ where

\smallskip

(i)  $h_l, h_r:\,B\to H$ are two group homomorphisms.

\smallskip

We also put $h^{\pm}:=h_l\pm h_r$ and denote the image
of $b\in B$ with respect to $h_l$ (resp. $h_r, h^{\pm})$ as $h_{b,l}$
(resp. $h_{b,r}, h^{\pm}_b)$.

\smallskip

(ii)  $\psi :\,B\to K^*$ is also a group homomorphism.

\smallskip

(iii)  $(\,, ):\,B\times B\to K^*$ is a symmetric pairing.

\smallskip

These data must satisfy the following condition: for all $b_1,b_2\in B$
$$
h_{b_2}^-(b_1)=(b_1,b_2)^2\alpha (h_{b_1,l},h_{b_2,l})\,
\alpha (h_{b_1,r},h_{b_2,r})^{-1}.
\eqno(2.12)
$$
\endproclaim

\medskip

{\bf 2.3.1. Remark.} Moduli space of quotients $T(H,\alpha )/B$
locally splits into a product of the classical moduli space of
commutative tori
$T(H,1 )/B$ and the space of quantization parameters $\alpha$
(which in a sense also are ``hidden periods'': cf. our discussion
of Connes' treatment of bad equivalence relations).

\smallskip
When $K=\bold{C}$, existence of sufficiently many theta functions
is equivalent to the algebraicity of $T(H,1 )/B$ which becomes an
abelian manifold. Multipliers of such theta functions
satisfy Riemann symmetry and positivity conditions.

\smallskip
Relations (2.12) represent an extension of the
symmetry conditions to our enlarged moduli space.
For a quantum version of positivity conditions, see Theorem 2.6.1 b)
below.
 
\medskip

{\bf 2.4. Automorphy factors.} For any theta multiplier
$\Cal{L}$ and period $b\in B$, {\it the automorphy factor}
$j_{\Cal{L}}(b)$ is, by definition, the following linear endomorphism
of any of the function spaces $Al, Af, An$ of $T(H,\alpha )$:
$$
j_{\Cal{L}}(b):\ \Phi\mapsto
\psi (b)\,(b,b)\,e(h_{b,l})\,\Phi\,e(h_{b,r})^{-1}.
\eqno(2.13)
$$
Clearly, it is invertible.

\medskip

\proclaim{\quad 2.5. Proposition} For $\Cal{L}$ fixed,
the map $b\mapsto j_{\Cal{L}}(b)^{-1}\circ b^*$
is a group homomorphism.

\smallskip

It is injective if $\roman{Ker}\,h^{-}= 0.$

\endproclaim

\smallskip

{\bf Proof.} We must check that
$$
j_{\Cal{L}}(b_1+b_2)^{-1}\circ (b_1+b_2)^*=
j_{\Cal{L}}(b_1)^{-1}\circ b_1^*\circ j_{\Cal{L}}(b_2)^{-1}\circ b_2^*.
\eqno(2.14)
$$
We have in view of (2.13), (2.4), (2.5):
$$
j_{\Cal{L}}(b)^{-1}(e(h))=\psi (b)^{-1}\,(b,b)^{-1}\,
\varepsilon(h_{b,l})\alpha (h,h_b^+)\alpha (h_{b,r},h_{b,l})\,
e(h -h^{-}_b).
\eqno(2.15)
$$
Now apply both sides of (2.14) to the arbitrary $e(h)$
and then compare them using (2.12).
A somewhat lengthy but straightforward calculation gives (2.14).
 
\smallskip

When $h^{-}$ is injective, $h^{-}_b\ne 0$ for non--zero $b$,
so that
$j_{\Cal{L}}(b)^{-1}\circ b^* (e(h))\ne e(h).$ 

\medskip

{\bf 2.6. Theta functions and theta types.} {\it A (quantized) theta function 
with multiplier $\Cal{L}$} is a formal series $\theta \in Af\,(T(H,\alpha ))$
invariant with respect to the transformation group
$$
\{j_{\Cal{L}}(b)^{-1}\circ b^*\,|\,b\in B\}.
$$ 
In other words,
$\theta$ must satisfy the functional equations
$$
b^*(\theta )=\psi (b)\,(b,b)\,e(h_{b,l})\,\theta\,e(h_{b,r})^{-1}
\eqno(2.16)
$$
for all $b\in B.$ Clearly, theta functions with multiplier $\Cal{L}$
form a linear space which we denote $\Gamma (\Cal{L})$. This notation
is supposed to remind the case of usual abelian
varieties where we deal with invertible sheaves and their sections.

\smallskip

Actually, different multipliers may produce the same space of
theta functions or even homomorphisms $b\mapsto j_{\Cal{L}}(b)^{-1}\circ b^*$.
Consider, for example, the case of
a commutative torus, where $\alpha\equiv 1.$ Then $\Gamma (\Cal{L})$
depends on $h_l,h_r$ only via their difference $h^-=h_l-h_r.$
Moreover, if $(\,,)^{\prime}$ is another symmetric pairing
such that $\varphi (b_1,b_2):=(b_1,b_2)^{\prime}(b_1,b_2)^{-1}$ takes values in $\{\pm 1\}$,
then $\varphi (b,b)$ is multiplicative in $b$, and we may replace
$(\psi ,(\,,))$ by $(\psi^{\prime} ,(\,,)^{\prime})$ where
$\psi^{\prime}(b)=\phi (b,b) \psi (b).$

\smallskip

Generally, from (2.15) one sees that if
$j_{\Cal{L}^{\prime}}(b)^{-1}\circ b^*=j_{\Cal{L}}(b)^{-1}\circ b^*$
for all $b$, then $h^-=h^{\prime -}$ and moreover,
$$
(b,b)^{-1}\alpha (h_{b,r},h_{b,l}) =\pm 
(b,b)^{\prime -1}\alpha (h_{b,r}^{\prime},h_{b,l}^{\prime}),
$$
$$
\psi (b)^{-1}\alpha (h,h^+_b) = \pm \psi^{\prime}(b)^{-1}\alpha 
(h,h^{\prime+}_b) .
$$

\smallskip

We will call two multipliers equivalent, if they have the same
space of theta--functions. (This definition is reasonable only
for ample multipliers, see below).
An equivalence class of multipliers $L$ will be called
{\it a theta type}. The space $\Gamma (\Cal{L})$ depends only
on this class and can be denoted also $\Gamma (L).$

\medskip

\proclaim{\quad 2.6.1. Theorem} a) We have 
$$
\roman{dim}\,\Gamma (\Cal{L}) = [H:h^-(B)].
\eqno(2.17)
$$
b) Assume that $K$ is a normed field.
Then all theta functions
of type $\Cal{L}$ are analytic if $[H:h^-(B)] <\infty$
and 
$$
\roman{log}\,|(b,b)\,\alpha (h_{b,l},-h_{b,r})|
\eqno(2.18)
$$ 
is a positively defined
quadratic form on $B$. 

\smallskip

In particular, assume that $B$ is free and
the quantization parameter $\alpha$ is unitary,
$|\alpha (h,h^{\prime})| \equiv 1.$ Then
this condition means that $\roman{rk}\,B = \roman{rk}\,H$,
$B$ is discrete in $T(H,1)(K)$ and $\roman{log}\,|(b,b)|$
is positively defined.
\endproclaim 

\smallskip

{\bf Proof.} Let $\theta =\sum_{h\in H} a_he(h),\,a_h\in K,$ $b\in B.$
We have
$$
b^*(\theta )=\sum_{h\in H}a_hh(b)e(h)=\sum_{h\in H}
a_{h+h^-_{b}}\,(h+h^-_{b})(b)\,e(h+h^-_{b}),
\eqno(2.19)
$$
whereas the right hand side of (2.16) is:
$$
\psi (b) (b,b)\sum_{h\in H}a_h \,\varepsilon (h_{b,r})\,
\alpha (h_{b,r},h_{b,i})\, \alpha (h_{b}^+,h)\, e(h+h^-_b).
\eqno(2.20)
$$
In (2.19) we replace $h^-_b(b)$ by 
$(b,b)^2\varepsilon (h_{b,l}) \varepsilon (h_{b,r})$
(see (2.12)). Comparing coefficients of (2.19)
and (2.20), we see that (2.16) is equivalent to
$$
a_{h+h^-_b}=a_h\, \psi (b)\,h(b)^{-1} (b,b)^{-1}
\varepsilon (h_{b,l})\,
\alpha (h_{b,r},h_{b,l})\, \alpha (h_{b}^+,h) 
\eqno(2.21)
$$
for all $h\in H$ and $b\in B.$
Thus one can arbitrarily choose values $a_h$ for all $h$ in
a system of representatives of $H/h^-(B)$ and then uniquely
reconstruct $\theta$. This proves the first statement of the theorem.
It also shows that if $\Gamma (\Cal{L})$ is not finite dimensional,
it necessarily contains non--analytic functions.
\smallskip

Assume now that $\Gamma (\Cal{L})$ is finite dimensional.
Then on each coset $h+h^-(B)$ we have
$$
\roman{log}\,|a_{h+h^-_b}| =
\roman{log}\,|a_h| -  \roman{log}\,|(b,b)\,\alpha (h_{b,r},h_{b,l})|
+ \roman{log}\,|\psi (b)\, h(b)^{-1}\alpha (h_{b}^+,h)|.
$$
The second summand in the right hand side is quadratic in $b$
whereas the third is linear. Hence analyticity follows from
the positive definiteness of (2.18).

\medskip

{\bf 2.6.2. Multiplication of theta functions.} Generally, we can multiply
analytic functions, but not formal ones. We will call
a theta multiplier $\Cal{L}$ {\it analytic}, if $\Gamma (\Cal{L})$
consists of analytic functions, and {\it ample}, if
$\Cal{L}$ satisfies conditions of Theorem 2.6.1 b).

\smallskip

We will call two analytic theta multipliers $\Cal{L}_i= (h_l^{(i)},h_r^{(i)},\psi_i , (\,,)_i) ,\,
i=1,2,$ {\it composable} (in this order) if
$h_l^{(2)}=h_r^{(1)}.$ Define their {\it product} as
$$
\Cal{L}_1\otimes\Cal{L}_2 =\Cal{L}:=
(h_l^{(1)},\,h_r^{(2)},\,\psi_1\psi_2,\, (\,,)_1(\,,)_2) .
$$
A straightforward calculation shows that if $\Cal{L}_i$ are ample,
$\Cal{L}$
is ample as well, and the product of theta
functions produces a well defined map
$$
\Gamma (\Cal{L}_1)\otimes \Gamma (\Cal{L}_2)\to
\Gamma (\Cal{L}):\ \theta_1\otimes\theta_2\mapsto \theta_1\theta_2.
\eqno(2.22)
$$
\smallskip

One can call two ample theta types $L_1$ composable, if
they contain pairs of composable multipliers. The multiplication in (2.22) 
does not
depend on the choice of such a pair, but I did not check that
the product type $L$ cannot change.

\smallskip

{\bf 2.6.3. Quantized abelian varieties.} Assume that $(T(H,\alpha ),B)$ admits ample theta multipliers. Then we consider $\oplus_L \Gamma (L)$
where $L$ runs over all theta types with $\psi (b)\in \{\pm 1\}$, together with
partial multiplication defined above, as a quantized version of the
graded coordinate ring of an abelian variety.
In the classical case, it is graded by the symmetric
elements of $\roman{Pic}$ lying in the effective cone (see 2.7.2 below).

\medskip

{\bf 2.7. Functorial properties of theta functions.} Consider
a morphism of non--commutative tori $F:\,T(H_2,\alpha_2 )\to T(H_1,\alpha_1)$.
As in Proposition 2.1.1 , let $F^*(e(h))=a_he(f(h)),\,a_h\in K^*$,
$f: H_1\to H_2.$ The map $f$ induces a morphism of commutative tori
which we also denote $F$:
$T(H_2,1)\to T(H_1,1)$.   
Let
$B_i\subset T(H_i,1)(K)$ be two period lattices such that
$F(B_2)\subset B_1.$

\smallskip

Choose a theta multiplier 
$\Cal{L}_1 = (h_l,h_r,\psi , (\,,))$ for $T(H_1,\alpha_1)$
and $B_1.$ We will show, how to produce from it
a new theta multiplier
$$
\Cal{L}_2=F^*(\Cal{L}_1) = (h_l^{\prime},h_r^{\prime},
\psi^{\prime}, (\,,)^{\prime})
$$
for $T(H_2,\alpha_2), B_2,$ 
if the following condition holds:
\smallskip

 {\it the characteristic of $F$ is 1, that is,
$a_{h+g}=a_ha_g$ (cf. end of the proof of Proposition 2.2.1),
and $f$ is compatible with $\alpha_1$ and $\alpha_2$,
and not just their squares.}

\smallskip

\smallskip

The map $f: H_1\to H_2$ induces a map $F:\,\roman{Hom}\,(H_2,K^*)\to
\roman{Hom}\,(H_1,K^*)$ whose restriction on $B_2$ sends it to $B_1$.
Put for $b, b_1, b_2\in B_2$:
$$
h_{l,r}^{\prime}=f\circ h_{l,r}\circ F :\ B_2\to H_2,
\eqno(2.23)
$$
so that $h_{b,l}^{\prime}=f(h_{F(b),l})$ and similarly for 
$h_r^{\prime}, h^{\prime \pm},$
$$
\psi^{\prime}(b)=\psi (F(b))\,a_{h^-_{F(b)}},
\eqno(2.24)
$$
$$
(b_1,b_2)^{\prime}=(F(b_1),F(b_2)) .
\eqno(2.25)
$$
\proclaim{\quad 2.7.1. Theorem} The data (2.23)--(2.25)
constitute a theta multiplier for $T(H_2,\alpha_2)$, $B_2$
such that $F^*(\Gamma (\Cal{L}))\subset \Gamma (F^*(\Cal{L})).$
\endproclaim

{\bf Proof.} We have first to check 
that  (2.12) holds for $F^*(\Cal{L})$.
For $b_1,b_2\in B_2$, the left hand side becomes
$$
h^{\prime -}_{b_2}(b_1)=f(h^-_{F(b_2)})(b_1)=h^-_{F(b_2)}(F(b_1)).
$$
Furthermore, if the characteristic of $F$ is 1,
the right hand side can be rewritten as
$$
(b_1,b_2)^{\prime 2}\alpha_2(f(h_{F(b_1),l}),f(h_{F(b_2),l}))
\alpha_2(f(h_{F(b_1),r}),f(h_{F(b_2),r}))^{-1}
$$
$$
=(F(b_1),F(b_2))^2\alpha_1(h_{F(b_1),l},h_{F(b_2),l})
\alpha_1(h_{F(b_1),r},h_{F(b_2),r})^{-1}.
$$
Applying (2.12) to $F(b_1),F(b_2)$ in lieu of $b_1,b_2$, we see that
both expressions coincide.

\smallskip

Let us now check that  $F^*(\Gamma (\Cal{L}))\subset \Gamma (F^*(\Cal{L})).$

\smallskip

Choose a formal function
$$
\theta =\sum_{h\in H_1} c_he(h),\ e=e_{H_1,\alpha_1}.
$$
According to (2.21), it belongs to $\Gamma (\Cal{L})$
iff the following conditions are satisfied
for all $h\in H_1,\,b\in B_1$: 
$$
c_{h+h^-_b}=c_h\, \psi (b)\,h(b)^{-1} (b,b)^{-1}
\varepsilon_1(h_{b,l})\,
\alpha_1(h_{b,r},h_{b,l})\, \alpha_1(h_{b}^+,h). 
\eqno(2.26)
$$
Notice that our notation slightly differs from (2.21):
$a_h$ is now reserved for $F^*(e(h))=a_he(f(h))$
as in the first lines of 2.7. 

\smallskip

Put now 
$$
F^*(\theta ) = \sum_{g\in H_2} C_ge(g),\quad e=e_{H_2,\alpha_2} .
$$
We have $C_g=0$, if $g\notin f(H_1)$; otherwise
$$
C_g=\sum_{h\in f^{-1}(g)}c_ha_h.
\eqno(2.27)
$$
To prove that $F^*(\theta )\in \Gamma (F^*(\Cal{L}))$,
we will check that $C_g$ satisfy analogs of 
relations (2.26) written for all $g\in H_2$
and all $b\in B_2.$

\smallskip

Consider first the case $g\notin f(H_1).$ Then
the relevant analog of (2.26) says that
$C_{g+h^{\prime -}_b}$ must be proportional to $C_g$
that is, zero. This is indeed true,
because in view of (2.23),
$h^{\prime -}_b=f(h^-_{F(b)})$ and thus $g+h^{\prime -}_b\notin f(H_1).$

\smallskip

Now assume that $g\in f(H_1)$, fix also $b\in B_2$, and
write down separately both sides of the relevant
case of (2.26). Since characteristic of $F$ is 1, the left hand side
can be written as
$$
C_{g+h^{\prime -}_b}=\sum_{h\in f^{-1}(g)}
c_{h+h^-_{F(b)}}a_h\,a_{h^-_{F(b)}}
$$
$$
=\sum_{h\in f^{-1}(g)} c_ha_h\, \psi (F(b)) a_{h^-_{F(b)}}\,
h(F(b))^{-1} (F(b),F(b))^{-1}\,
$$
$$
\times
\varepsilon_1(h_{F(b),l})\,
\alpha_1(h_{F(b),r},h_{F(b),l})\,\alpha_1(h^+_{F(b)},h).
\eqno(2.28)
$$
Here we have rewritten $c_{h+h^-_{F(b)}}$ using (2.26).

\smallskip

The right hand side is
$$
\left(\sum_{h\in f^{-1}(g)}c_ha_h\right)\psi^{\prime}(b)\, g(b)^{-1}
(F(b),F(b))^{-1}\,
 \varepsilon_2(h^{\prime -}_{b,l})\,
\alpha_2(h^{\prime -}_{b,r}, h^{\prime -}_{b,l})\,
\alpha_2(h^{\prime +}_{b},g).
\eqno(2.29)
$$
We can now compare (2.28) and (2.29) term by term
using (2.23)--(2.25), and convince ourselves
that they coincide. In particular, we
use the identities $g(b)=h(F(b))$
and $g=f(h)$ for any $h\in f^{-1}(g).$

\medskip

{\bf 2.7.2. Examples.} All morphisms of tori described in 2.1.3
have characteristic 1. Therefore, complementing
these tori by compatible period lattices, we obtain
quantized versions of many standard morphisms of abelian varieties.

\smallskip

In particular, $[-1]^*$ acts on theta multipliers
for $T(H,\alpha ),B$ by simply inverting $\psi$.
Hence symmetric theta multipliers correspond
to $\psi : B\to \{\pm 1\}.$

\smallskip

Mumford's morphism for abelian varieties, induced by
the toric morphism (2.11), is the starting point
for his study of homogeneous coordinate rings
of abelian varieties. Our preparatory work
allows us to hope that at least part of this study can be
extended to the quantum case.
\bigskip

\bigskip

\centerline{\bf Appendix. Commutative Geometry as Noncommutative Geometry}

\medskip

The main goal of this Appendix is to illustrate
Connes approach to noncommutative geometry in an
algebraic geometric context and to make convincing
our claim that an elliptic curve ``is'' a
two--dimensional noncommutative torus. 

\smallskip

We take as our starting point Connes explanations
about how to treat as a non--commutative space
quotient of a ``commutative space'' by an equivalence
relation or equivalence groupoid: see [Co], II.2--II.5.
This viewpoint is complementary to the more popular
in algebraic geometry deformation paradigm.

\smallskip

First, a reminder about groupoids.

\smallskip

Let  $U$
be a set. Classically, an equivalence relation $\sim$ on 
$U$ is given by its graph $R\subset U\times U,\,R:=\{(a,b)\,|\,a\sim b\}$
which satisfies the three conditions:
\smallskip

{\it Reflexivity:} 
$$
a\sim a\ \Longleftrightarrow\ \Delta_U\subset R;
$$

\smallskip

{\it Symmetry:} 
$$
(a\sim b\,\Leftrightarrow b\sim a)\ \Longleftrightarrow\ 
s_{12}(R)=R;
$$
\smallskip

{\it Transitivity:}
$$
((a\sim b) \& (b\sim c) \Rightarrow a\sim c)\
\Longleftrightarrow\ \roman{pr}_{13}[(R\times U)\cap (U\times R)]\subset R.
$$
\smallskip

All of this can be rephrased as follows: there exists
a category with the set of objects $U,$ set of morphisms $R,$
such that $R\to U\times U$ is the map $f\mapsto$ (source of $f$,
target of $f$), and in addition, every morphism is
an isomorphism, and all automorphism groups are
trivial.

\smallskip

Consider now a diagram $R\to U\times U$ satisfying this
description with the last condition deleted so that the
automorphism groups can now be arbitrary. We will call
such a diagram {\it an equivalence groupoid}
(on the set $U$). Of course,
an equivalence groupoid $R\to U\times U$ comes together
with the identity map $U\to R:\, a\mapsto \roman{id}_a$
and the associative multiplication map $R\times_UR\to R$
satisfying the usual categorical axioms which reduce to
the reflexivity, symmetry, and transitivity for the usual
equivalence relations. Notice that the image of $R$
is in fact an equivalence relation, and the respective quotient is
the set of isomorphism classes of objects.

\smallskip

Thus the basic difference between equivalence groupoids and
equivalence relations on sets can be demonstrated on one--point
sets $U=\{*\}$: in this case $R$ is simply a group.
In the framework of homotopy theory, the respective
quotient object $\{*\}/R$ is represented by the
classifying space $BR$. Stacks provide a categorical
context for constructing such quotients in algebraic geometry.

\smallskip 

In fact, the notion of equivalence groupoid was formulated
in such a way that it readily generalizes
to the case when $R\to U\times U$ is a diagram in
an arbitrary category with products, e.g. schemes.
One should imagine $U$ as an atlas and $R$ as gluing
rules, so that the geometric object we are
interested in is symbolically $U/R$.

\smallskip

Connes prescription, roughly speaking, consists
in reducing geometric study of $X=U/R$ to the 
representation theory of a non--commutative
ring $\Cal{R}_X$. The heuristic rule for
constructing $\Cal{R}_X$ can be stated as follows.
Let us write elements (points) of $R$ as morphisms
$j:\,u\to u^{\prime}$ so that $R\to U\times U$
maps such a point to $(u,u^{\prime}).$ Then $\Cal{R}_X$
consists of certain functions $f$ on $R$
endowed with convolution multiplication:
$$
(f*g)(k: u\to u^{\prime\prime})=
\sum_{(i,j):\,ij=k}f(j: u\to u^{\prime})\,
g(i: u^{\prime}\to u^{\prime\prime})\,.
\eqno(A.1)
$$
Of course, correct choice of the class of functions and making sense of
the convolution multiplication may present a problem,
but in algebraic geometry it is clear what to start with at least when $U$
is affine.
\smallskip

Let us illustrate such a setup by several examples.

\smallskip

{\it Example 1.} If $X$ is the quotient $\{*\}/G$ where $G$
is a finite group, then $\Cal{R}_X$ is the group ring of $G$.
Representation theory of $\Cal{R}_X$ is essentially
$K$--theory of the classifying stack of $G$, in accordance
with the common wisdom.

\medskip

{\it Example 2.} If $X$ is the quotient of an affine scheme
$U=\roman{Spec}\,A$ by the action of a group $G$, then
$\Cal{R}_X$ should contain at least the twisted product
of $A$ with a version of the group ring of $G$.

\smallskip

When $U=\bold{G}_m$ with coordinate $y$ 
and $G=\bold{Z}$ with generator $x$ acting on $y$
as multiplication by $q$, then in this twisted
product we have $xyx^{-1}=qx$ that is, the basic
relation of the noncommutative torus.
On the other hand, $U/G$ makes sense as an elliptic
curve in analytic geometry, if $|q|\ne 1.$
 
\medskip

{\it Example 3.} A projective scheme $\roman{Proj}\,A$
is the quotient of its cone $\roman{Spec}\,A$ (with vertex
deleted) by the action of $\bold{G}_m$ determining
the grading. Studying representations of the
twisted product of $A$ with the group ring of $G_m$
is equivalent to studying graded $A$--modules.
Deleting the vertex boils down to taking
the quotient of the category of graded modules by
the subcategory of modules with finite number
of nonvanishing components. This is the familiar
Serre's picture which is archetypal in the following sense:
after finding an appropriate ring $\Cal{R}_X$,
one proceeds to establish an equivalence of categories
$\Cal{R}_X$--$\roman{Mod} \to \roman{Coh}_X$, or
its localized and/or derived version.

\bigskip

\centerline{\bf References}

\medskip

[AP] D.~Arinkin, A.~Polishchuk. {\it Fukaya category
and Fourier transform.} Preprint math.AG/9811023

\smallskip

[Bar] S.~Barannikov. {\it Extended moduli spaces and mirror
symmetry in dimensions $n>3.$} Preprint math.AG/9903124

\smallskip

[BarK] S.~Barannikov, M.~Kontsevich. {\it Frobenius manifolds
and formality of Lie algebras of polyvector fields.} Int. Math. Res. 
Notices, 4 (1998), 201--215.

\smallskip

[BEG] V.~Baranovsky, S.~Evens, V.~Ginzburg. {\it Representations
of quantum tori and double--affine Hecke algebras.} Preprint
math.RT/0005024

\smallskip

[Bat] V.~Batyrev. {\it Dual polyhedra and the mirror symmetry
for Calabi--Yau hypersurfaces in toric varieties.} Journ. Alg. Geom.,
3 (1994), 493--535.

\smallskip

[Co] A.~Connes. {\it Noncommutative geometry.} Academic Press,
1994.

\smallskip

[De] P.~Deligne. {\it Local behavior of Hodge structures at infinity}. 
In: Mirror Symmetry II,
ed.~by B.~Greene and S.~T.~Yau, AMS--International Press, 1996, 683--699.

\smallskip 

[Fu] K.~Fukaya. {\it Mirror symmetry of abelian variety and
multi theta functions.} Preprint, 1998.

\smallskip

[Giv] A.~Givental. {\it Equivariant Gromov--Witten invariants.}
Int. Math. Res. Notes, 13 (1996), 613--663.

\smallskip

[GolLO] V.~Golyshev, V.~Lunts, D.~Orlov. {\it Mirror symmetry for
abelian varieties.} Preprint math.AG/9812003

\smallskip

[Gr1] M.~Gross. {\it Special Lagranfian fibrations I: Topology.}
 alg-geom/9710006

\smallskip

[Gr2] M.~Gross. {\it Special Lagranfian fibrations II: Geometry.}
 math.AG/9809072

\smallskip

[Ko]  M.~Kontsevich. {\it Homological
algebra of Mirror Symmetry.} Proceedings of the ICM
(Z\"urich, 1994), vol. I, Birkh\"auser, 1995, 120--139.
Preprint alg-geom/9411018.



\smallskip

[LYZ] N.~C.~Leung, Sh.-T.~Yau, E.~Zaslow. {\it From special
Lagrangian to Hermitian--Yang--Mills via Fourier--Mukai
transform.} Preprint math.DG/0005118.

\smallskip

[Ma1] Yu.~Manin. {\it Quantized theta--functions.} In: Common
Trends in Mathematics and Quantum Field Theories (Kyoto, 1990), 
Progress of Theor. Phys. Supplement, 102 (1990), 219--228.

\smallskip

[Ma2] Yu.~Manin. {\it Frobenius manifolds, quantum cohomology, and moduli
spaces.} AMS Colloquium Publications, vol. 47, Providence, RI, 1999,
xiii+303 pp.

\smallskip

[Ma3] Yu.~Manin. {\it Moduli, Motives, Mirrors.} Plenary talk
at 3rd European Congress of Mathematicians, Barcelona, 2000,
Preprint.

\smallskip

[MirS1] S.--T.~Yau, ed. {\it Essays on Mirror Manifolds.} International Press
Co., Hong Cong, 1992.

\smallskip

[MirS2] B.~Greene, S.~T.~Yau, eds. {\it Mirror Symmetry II.},
AMS--International Press, 1996.

\smallskip

[Mu] D.~Mumford. {\it On the equations defining abelian
varieties I.} Inv. Math. 1 (1966), 355--374.

\smallskip

[Po1] A.~Polishchuk. {\it Massey and Fukaya products on elliptic
curve.} Preprint math.AG/9803017

\smallskip

[Po2] A.~Polishchuk. {\it Homological mirror symmetry with higher
products.}  Preprint math.AG/9901025



\smallskip

[PoZ] A.~Polishchuk, E.~Zaslow. {\it Categorical mirror
symmetry: the elliptic curve.} Adv. Theor. Math. Phys.,
2 (1998), 443--470. 
Preprint math.AG/9801119

\smallskip

[RiS] M.~A.~Rieffel, A.~Schwarz. {\it Morita equivalence
of multidimensional tori.} Preprint math.QA/9803057

\smallskip

[So] Y.~Soibelman. {\it Quantum tori, mirror symmetry
and deformation theory.} Preprint, 2000.

\smallskip

[StYZ] A.~Strominger, S.--T.Yau, E.~Zaslow. {\it Mirror symmetry
is $T$--duality.} Nucl.~Phys. B 479 (1996), 243--259.





\smallskip

[We] A.~Weinstein. {\it Classical theta functions and quantum tori.}
Publ. RIMS, Kyoto Univ., 30 (1994), 327--333.



\enddocument